\documentclass[11pt]{amsart}
\usepackage{geometry}                
\usepackage{graphicx}
\usepackage{amssymb}
\usepackage{epstopdf}
\usepackage{mathrsfs}
\title{Classifying Several Classes of Leibniz Algebras}
\author{Chelsie Batten Ray, Allison Hedges and Ernest Stitzinger}

\begin{document}
\maketitle

\vspace{5mm}

ABSTRACT

\vspace{5mm}

We extend results related to maximal subalgebras and ideals from Lie to Leibniz algebras. In particular, we classify minimal non-elementary Leibniz algebras and Leibniz algebras with a unique maximal ideal. In both cases, there are types of these algebras with no Lie algebra analogue. We also give a classification of E-Leibniz algebras which is very similiar to its Lie algebra counterpart. Note that a classification of elementary Leibniz algebras has been shown in [3].

\vspace{5mm}

I. PRELIMINARIES

\vspace{5mm}

Loday introduced Leibniz algebras as a noncommutative generalization of Lie algebras. Lie algebra results have been extended to these new algebras. Properties of the Frattini subalgebra and ideal are studied in [3] for Leibniz algebras. In particular, elementary Leibniz algebras, those algebras with the property that the Frattini ideal is 0 for every subalgebra, were classified over algebraically closed fields, extending Lie algebra results of Towers  [12], [13] and Towers and Varea [14]. Algebras closely related to elementary Lie algebras are E-Lie algebras and minimal non-elementary Lie algebras, which have have been investigated in [13]. We extend these results to Leibniz algebras. In the minimal non-elementary case, the classification contains many new algebras. Similar in definition to the Frattini subalgebra is the Jacobson radical, which is the intersection of all maximal ideals. This concept appears in [5] and [9] for Lie algebras. In particular, Lie algebras with a unique maximal ideal are classified in [5]. We investigate these ideas in Leibniz algebras, again getting new cases in the classification result.

\vspace{5mm}

Let A be an algebra. The intersection of all maximal subalgebras, F(A), is called the Frattini subalgebra and the largest ideal of A contained in F(A), $\Phi$(A), is called the Frattini ideal. For Lie algebras these concepts have been widely studied, see [13]. If A is a solvable Lie algebra or is over a field of characteristic 0, then F(A) is an ideal, but not generally. For Leibniz algebras F(A) is an ideal at characteristic 0, but not necessarily when A is solvable [3]. The intersection of all maximal ideals, J(A), is called the Jacobson radical. The nilradical, Nil(A), is the maximal nilpotent ideal of A; it exists by [6]. The sum of all minimal abelian ideals is denoted by Asoc(A). An algebra is elementary if $\Phi$(B)=0 for all subalgebras B of A and is minimal non-elementary if $\Phi$(B)=0 for all proper subalgebras of A but $\Phi$(A) $\neq 0$. A is called an E-algebra if $\Phi$(B)$ \subseteq \Phi$(A) for all subalgebras of B of A.

\vspace{5mm}

Following Barnes [1], we call an algebra L (left) Leibniz if left multiplication by each x
in L is a derivation. Many authors consider right Leibniz algebras instead. Thus an algebra L is Leibniz if x(yz)=(xy)z+y(xz) is an identity on L. The product need not be antisymmetric.

\vspace{5mm}

II. E-LEIBNIZ ALGEBRAS

\vspace{5mm}

L is an E-algebra if $\Phi$(B)$\subseteq \Phi$(L) for all subalgebras B of L. The following result is an extension of Proposition 2 of [11] to Leibniz algebra. The proof is the same as the original, hence we omit it. 

\vspace{5mm}

Theorem 1. Let L be a Leibniz algebra. Then L is an E-Leibniz algebra if and only if L/$\Phi$(L) is elementary.

\vspace{5mm}

Proposition 2. If L is solvable over a field of characteristic 0, then L is an E-Leibniz algebra.

\vspace{5mm}

Proof. Since L is solvable, L$^2$ is nilpotent. By Theorem 3.5 of [3] the result holds.

\vspace{5mm}

Lemma 3. If L=B$\bigoplus$C and B and C are elementary, then L is elementary.

\vspace{5mm}

Proof. Let S be a subalgebra of B$\bigoplus$C. We show that $\Phi$(S)=0. For x$\in$B+S,  x=b+c, where b$\in$B and c$\in$C and c=x-b $\in$C$\cap$(B+S)
The projection mapping from B+S onto (B+S)$\cap$ C, where x goes to c, has kernel B. Hence S/(B$\cap$S)$\cong$(B+S)/B$\cong$ (B+S)$\cap$C. Since C is elementary, $\Phi$((B+S)$\cap$C)=0.  Hence $\Phi$(S)$\subseteq$B$\cap$S. Similarly $\Phi$(S)$\subseteq$C$\cap$S. Hence $\Phi$(S)$\subseteq$B$\cap$C=0.

\vspace{5mm}

Theorem 4. Let L be a Leibniz algebra over K, an algebraically closed field of characteristic 0. Then L is an E- Leibniz algebra if and only if:

(1) L is solvable, or

(2)  L$\cong$ sl$_2$(K)$\bigoplus...\bigoplus$sl$_2$(K), or

(3) L=R+S, where R is a solvable ideal, S $\cong$ sl$_2$(K)$\bigoplus...\bigoplus$sl$_2$(K), and RS+SR is contained in $\Phi$(L).

\vspace{5mm}

Proof. If L is solvable, then L is an E-Leibniz algebra by Proposition 2. 
If L is as in (2), then L is a Lie algebra and L is elementary by Theorem 3.2 of [13].
Let L be as in (3). S is elementary as in (2). $\Phi$(L) is a nilpotent ideal in L, hence $\Phi$(L)
$\subseteq$ R. 
Thus, L/($\Phi$(L))$ \cong$R/$\Phi$(L)+S. Since RS+SR$\subseteq \Phi$(L), the previous sum is a direct sum. Thus, by Theorem 4.8 of [12], 0=$\Phi$(L/$\Phi$(L))=$\Phi$(S)$\bigoplus \Phi$(R/$\Phi$(L)). This shows that R/$\Phi$(L) is elementary using Proposition 2. Then L/$\Phi$(L) is elementary by Lemma 3 and L is an E-algebra by Theorem 1.

\vspace{5mm}

Conversely, suppose that L is an E-algebra. Then L/$\Phi$(L) is elementary from Theorem 1. L/$\Phi$(L) is the direct sum of its radical, R, and a semisimple ideal S  $\cong$ sl$_2$(K)$\bigoplus...\bigoplus$sl$_2$(K) by Theorem 4.3 of [3], either of which may be 0. If S=0, then L is solvable and (1) holds. If R=0, the L is the direct sum of copies of sl$_2$(K) as in (2). If neither R nor S is 0, then RS+SR $\subseteq \Phi$(L) as a consequence of Theorem 4.3 of [3].

\vspace{5mm}

The following result addresses the special case in which L is a perfect Leibniz algebra.

\vspace{5mm}

Corollary 5. Let L be a perfect Leibniz algebra (L$^2$=L) over K, an algebraically closed field of characteristic 0. Then L is an E-Leibniz algebra if and only if  L $\cong$ sl$_2$(K)$\bigoplus...\bigoplus$sl$_2$(K).

\vspace{5mm}

Proof. Let L=R+S be the Levi decomposition for L as in [2], where R is the radical and S is a semisimple subalgebra. Then L$^2$=S$^2$+RS+SR+R$^2$. R=RS+SR+R$^2$ since L is perfect. Since L is an E-algebra, RS+SR $\subseteq \Phi$(L). Hence R$^2$+$\Phi$(L)=R and R$^2$+$\Phi$(L)+S=L. Hence R$^2$+S=L since no proper subalgebra can supplement $\Phi$(L). Thus R$^2$=R, which implies that R=0 since R is solvable. Hence L=S and L is Lie. Then, by Corollary 4.5 of [13], L$\cong$ sl$_2$(K)$\bigoplus...\bigoplus$sl$_2$(K).

\vspace{5mm}

III. MINIMAL NON-ELEMENTARY LEIBNIZ ALGEBRAS

\vspace{5mm}

An algebra, L, is called minimal non-elementary if L is not elementary but all proper subalgebras of L are elementary.  In [5], conditions for a Lie algebra to be minimal non-elementary are found when L$^2$ is nilpotent. The next result is the Leibniz algebra version. Note that cases 3 and 4 have no Lie algebra counterpart and case 1 has a case not found in Lie algebras.   When L$^2$ is nilpotent, $\Phi$(L)=0 if and only if L$^2 \subseteq$ Asoc(L) and L$^2$ is complemented in L by Proposition 3.1 of [3]. We often use this fact in the following proof.

\vspace{5mm}
 Theorem 6. Let L be a finite dimensional Leibniz algebra over an algebraically closed field. Suppose that L$^2$ is nilpotent. L is minimal non-elementary if and only if: 
 \vspace{2mm}

(1) L is three dimensional non-nilpotent with basis {x,y,z} and non-zero multiplication as:

\vspace{2mm}

(a) xz=cz, xy=cy+z, zx=0 and yx=0, where c is a non-zero scalar, or

(b) xz=cz, xy=cy+z, zx=-cz, yx=-cy-z where c is a non-zero scalar, or  

\vspace{2mm}

(2) L is Heisenberg, or

\vspace{2mm}

(3) L is generated by a, where a$^2 \neq 0$, L is nilpotent and dim L  $\geq 2$, or

\vspace{5mm}

(4) L is the four dimensional non-nilpotent algebra generated by a,b,x and y with multiplication:

ax=$\alpha$ x, ay=$\alpha$y, bx=$\beta$x, by=$\beta$y, xa=-$\alpha$x, ya=0, xb=y-$\beta$x and yb=0 where $\alpha, \beta $  are non zero scalars.

\vspace{5mm}

Proof.  Let L  be minimal non-elementary. L is solvable by the conditions. Suppose that L is not nilpotent. Then $\Phi$(L)$\neq$L$^2$ and there exists a maximal subalgebra, M, of L such that L=L$^2$+M. Let B be an algebra of minimum dimension such that L=-L$^2$+B. By Lemma 7.1 of [12], L$^2 \cap$B $\subseteq \Phi$(B)=0. L$^2$ is nilpotent and elementary, hence it is abelian. Clearly B is also abelian.

\vspace{5mm}

Suppose that dim B $> $1. For any a$\in$B , let H(a)=L$^2$+(a). Then $\Phi$(H(a))=0 since all proper subalgebras are elementary. Then L$^2\subseteq$ Nil(H(a))=Asoc(H(a)). L$^2$ is abelian since it is elementary and nilpotent, and it is completely reducible under the action of a. On each minimal ideal, either R$_a$=-L$_a$ or R$_a$=0 [1]. Hence the minimal ideals are one dimensional  eigenspaces for L$_a$ and R$_a$ acting on L$^2$ and L$_a$ and R$_a$ are simultaneously diagonalizable on L$^2$.  This holds for all a in B, and since the left multiplications, L$_a$, commute, they are simultaneously diagonalizable. 

\vspace{5mm}

Consider B acting on W=L$^2$ by left multiplication and decompose L$^2$ as the direct sum of weight modules, \{x $\in$ W: ax=$\alpha_a$x for all a in B\}  Each of these weight modules is invariant under R$_b$, b $\in$B. 
If there is more than one weight module, then, for each weight $\alpha$, let B$_\alpha$=W$_\alpha$+B. Since $\Phi($B$_\alpha$)=0, W$_\alpha \subseteq$ Asoc(B$_\alpha$)=Nil(B$_\alpha$) and W$_\alpha$ is completely reducible as a B-bimodule. Hence W $\subseteq$ Asoc(L). Since B complements W, Asoc(L) is complemented and $\Phi$(L)=0, a contradiction. 

Hence there is one weight module and each left multipliction is a scalar. Pick a$\in$ B and suppose that ax=$\alpha$x on W. Let W$_0$= \{x $\in$ W: xa=0\} and W$_1$=\{x$\in $W: xa= -$\alpha$x\}. Since W is the direct sum of  one dimensional a-invariant submodules, W=W$_0$+W$_1$. If, for each a $\in$ B, W=W$_0$ or W=W$_1$, then each right multiplication is a scalar on L$^2$, and L$^2\subseteq $Asoc(L). Again Asoc(L) is complemented in L and $\Phi$(L)=0, a contradiction. 
Thus assume there is an a $\in$ B such that neither W$_0$ nor W$_1$ is 0. W$_0$ is a submodule. If W$_1$ is a submodule, then by induction both components are completely reducible under B, and W $\subseteq$ Asoc(L). Thus Asoc(L) is complemented and $\Phi$(L) = 0, a contradiction. Hence assume there is an x in W$_1$ and  b in B such that xb is not in W$_1$. Let N=(x,xb), bx=$\beta$x, and set y=$\beta$x+xb. The following multiplications hold:
ax=$\alpha$x, ay=$\alpha$y, bx=$\beta$x, by=$\beta$y, xa=-$\alpha$x, ya=($\beta$x+xb)a=$\beta$xa-b(xa)=-$\alpha \beta$x+$\alpha \beta$x=0, xb=y-$\beta$x, and yb=($\beta$x+xb)b=$\beta$xb-b(xb) =$\beta$xb-$\beta$xb=0. Let C = (a,b).  Hence N is a C-bimodule and (y) is a submodule which is not complemented in N, for suppose that $\sigma$x+$\tau$y is C-invariant where $\sigma\neq 0$. Then ($\sigma$x+$\tau$y)a=-$\sigma \alpha$x which yields that $\tau=0$. Then $\sigma$x b= $\sigma$y-$\sigma \beta$ x and $\sigma$=0 since xb is not in W$_1$. Thus no complement exists. Hence Asoc(N+C) is not complemented in N+C and $\Phi$(N+C) $\neq $0. Thus B=C, N=W=L$^2$, dim B=2 and dim L$^2$=2 and the multipliction for B acting on L$^2$ is the one given in this paragraph.

\vspace{5mm}

\vspace{5mm}

Let B=(a,b) and L$^2$=(x,y). From the forgoing ax=$\alpha$x, ay=$\alpha$y bx=$\beta$x, by=$\beta$y, xa=-$\alpha$x, ya=0, yb=0, and xb=y-$\beta$x. N=(y) is a one dimensional submodule of L, and L$^2$ is not completely reducible under the action of L. This is the algebra in (4).

\vspace{5mm}

Suppose that dim B =1. Hence L=L$^2$+B and L$^2$ is abelian. We now show that there exists a chain of ideals 0 $\subseteq$ L$_1$...L$_{n-1} \subseteq$=L$^2 \subseteq$ L$_n$=L where dim L$_i$=i. If P $\subseteq$ Q$\subseteq$L$^2$ are ideals of L with Q/P irreducible under the action of L, then dim Q/P=1 since the action of x on Q/P determines the action of L on Q/P, and either R$_x$=-L$_x$ or R$_x$=0. Then any eigenvector of L$_x$ in Q/P must span Q/P and dim Q/P=1. Hence there exists a flag from 0 to L as claimed. Now M=L$_{n-2}$+B is a maximal subalgebra of L, hence $\Phi$(M) = 0 by assumption. Hence L$_{n-2} \subseteq$ Asoc(M) = Nil(M). Since L$^2$ is abelian and B is one dimensional, Asoc(M) $\subseteq$ Asoc(L). Thus L$_{n-2} \subseteq$ Asoc(L). If L$_{n-2} \neq $ Asoc(L), then Asoc(L)=L$_{n-1}$=Nil(L), and L splits over Nil(L)=Asoc(L); hence, $\Phi$(L)=0 from Theorem 3.1 of [3], and L is elementary. Otherwise, Asoc(L)=L$_{n-2}$ and Asoc(L) has co-dimension two in L. If Asoc(L) is not contained in $\Phi$(L), then L splits over Asoc(L) by Theorem 7.1 of [12]. Again $\Phi$(L)=0. and L is elementary. Thus assume that Asoc(L)=$\Phi$(L). Since all minimal ideals are one dimensional, they are eigenspaces for L$_x$ where B=$<$x$>$. Note that R$_x$=0 or R$_x$=-L$_x$ on each of these minimal ideals. Since L$^2$/Asoc(L) is one-dimensional, there exists a scalar $\alpha$ with (L$_x-\alpha I)^2$=0 and L$_x$-$\alpha I \neq 0$. Hence there exist y,z $\in$ L$^2$ with the following possible multiplications.

\vspace{5mm}

Case 1: Let xz=cz, xy=cy +z, zx=0, yx=-dcy+ez, where d=0 or 1 and c is a non-zero scalar since L is not nilpotent. From x(yx)=(xy)x+y(xx) we obtain that -cd +ce=ce. Hence d=0. From y(xx)=(yx)x+x(yx), we obtain that -cde-cd+ce=0. Hence e=0, and yx=0.

\vspace{5mm}

Case 2: Let xz=cz, xy=cy+z, zx=-cz, yx=-dcy+ez where d=0 or 1 and c $\neq$ 0. Using x(yx)=(xy)x +y(xx), we obtain that -dc=-c. Hence d=1. From y(xx)=(yx)x+x(yx), we obtain that cde=-cd. Hence e=-1, and yx=-cy-z.

\vspace{5mm}

Hence H=$<x,y,z>$ has $\Phi$(H)=$<z>$ and L=H, which is the algebra in (1).

\vspace{5mm}

Let L be nilpotent with all proper subalgebras elementary. If there exists an a $\in$ L with a$^2 \neq 0$, then the subalgebra B with basis {a, a$^2$,...,a$^n$} and aa$^n$=0 has $\Phi$(B)=B$^2$. Since a$^2 \in \Phi$(B), B=L, which is the algebra in (3). If no such a exists, then L is Lie and hence Heisenberg by Theorem 4.7 of [8], yielding case (2).

\vspace{5mm}

Conversely, in cases (1) and (2) and (4), clearly all proper subalgebras are elementary. Suppose that L is as in case (3). Then L=$<$a,a$^2$,...,a$^n>$ with aa$^n$=0. Since L is nilpotent, L$^2$=$\Phi(L)$. Now b= c$_1$a+c$_2$a$^2$+...+c$_n$a$^n$ is in $\Phi$(L) if and only if c$_1$=0. If the subalgebra b contains an element that is not in $\Phi$(L), then B+$\Phi$(L)=L since L/$\Phi$(L) is one dimensional. Hence B=L. Therefore all proper subalgebras of L are contained in L$^{2}$, which is abelian, and hence are elementary. Thus L satisfies the conditions of the theorem.

\vspace{5mm}

We turn to the Leibniz algebra version of a Lie algebra result of Towers [13].

\vspace{5mm}

Theorem 7. Let L be a Leibniz algebra over K, an algebraically closed field of characteristic 0. L is minimal non-elementary if and only if:

\vspace{2mm}

(1) L is three dimensional non-nilpotent with basis {x,y,z} and non-zero multiplication as:

\vspace{2mm}

(a) xz=cz, xy=cy+z, zx=0 and yx=0, where c is a non-zero scalar, or

(b) xz=cz, xy=cy+z, zx=-cz, yx=-cy-z where c is a non-zero scalar, or  

\vspace{2mm}

(2) L is Heisenberg, or

\vspace{2mm}

(3) L is generated by a, where a$^2 \neq 0$, L is nilpotent and dim L  $\geq 2$, or

\vspace{2mm}

(4) L is the four dimensional non-nilpotent algebra generated by a,b,x and y with multiplcation.
ax=$\alpha$ x, ay=$\alpha$y, bx=$\beta$x, by=$\beta$y, xa=-$\alpha$x, ya=0, xb=y-$\beta$x and yb=0 where $\alpha$ and $\beta $  are non zero scalars.

\vspace{5mm}

Proof: Suppose that L is minimal non-elementary. Then L is an E-algebra.

\vspace{5mm}

Suppose that L is not solvable. Hence there exists a k such that L$^{(k)}$=L$^{(k+1)}$ If k=1, then L is a perfect Lie algebra and L is elementary by Corollary 5, a contradiction. If k$\geq 2$,  then L$^{k}$ is perfect and L$^{(k)}\cong$ sl$_2$(K)$\bigoplus...\bigoplus$sl$_2$(K). If R is the radical of L, then, then L=R$\bigoplus$L$^{(k)}$. Since both summands are elementary, L is elementary by Lemma 3 , a contradiction. Hence L must be solvable. The result now follows from Theorem 6.

\vspace{5mm} 

IV. THE JACOBSON RADICAL

\vspace{5mm}

The Jacobson radical J(L), is the intersection of maximal ideals of L. This concept was considered in [5] and [9] when L is a Lie algebra. If L is nilpotent, then J(L)=$\Phi$(L), since then all maximal subalgebras are ideals [1]. Clearly J(L)$\subseteq$L$^2,$ since if x is not in L$^2$, then we can find a complementary subspace, M, of x in L that contains L$^2$ and, since L$^2\subseteq$ M, M is a maximal ideal of L and x is not in M.

\vspace{5mm}
If L is a linear Lie algebra, let R=Rad(L), and let Rad(L$^*$) be the radical of the associative envelope, L$^*$ of L. Then, by corollary 2 p. 45 of [7], L$\cap$ Rad(L$^*$) = all nilpotent elements of R and [R,L] $\subseteq$Rad(L$^*$).

\vspace{5mm}
Theorem 8. Let L be a Leibniz algebra and R=Rad(L) be the radical of L. Then LR+RL$\subseteq$N=Nil(L).

\vspace{5mm}

Proof: Let $\mathscr{L}$(L)=\{L$_x$: x$\in$N\}. The map $\pi$:L$\rightarrow \mathscr{L}$(L) is a homomorphism and $\mathscr{L}$(L) is a Lie algebra under commutation.. Also $\pi:$R$\rightarrow$R($\mathscr{L}$(L)), the radical of $\mathscr{L}$(L). By the result in the first paragraph, [$\mathscr{L}$(L),R($\mathscr{L}$(L))]$\subseteq $R($\mathscr{L}$(L)$^*$). Hence there exists an n such that [L$_{s_n}$,L$_{t_n}]$...[L$_{s_2}$,L$_{t_2}]$[L$_{s_1}$,L$_{t_1}]$=0, where s$_i \in$ L , t$_i$ R or s$_i \in$ R and t$_i\in$ L. Hence L$_{s_nt_n}$...L$_{s_1t_1}=0$. Hence $s_nt_n(...(s_1t_1(x))...)=0$ for all x$\in$L. Hence $s_nt_n(...(s_1t_1(s_0t_0))...)=0$ Therefore(LR+RL)$^{n+1}$=0. Since LR+RL is an ideal in L, LR+RL  $\subseteq$ N.

\vspace{5mm}

Proposition 9. If L is of characteristic 0, then L$^2\cap$R=LR+RL $\subseteq$N.

\vspace{5mm}

Proof: L=R+S as in the Levi decomposition for Leibniz algebras ([2]). Then L$^2=$S$^2+$LR+RL. and L$^2\cap$R=S$^2\cap$R+(LR+RL)$\cap$R=LR+RL since S$\cap$R=0 and LR+RL$\subseteq$R.

\vspace{5mm}

Lemma 10. If L is solvable, then J(L)=L$^2$.

\vspace{5mm}

Proof: If M is a maximal ideal of L, then L/M is abelian and, in fact, one-dimensional since L is solvable. Hence L$^2 \subseteq$ M for all M and L$^2 \subseteq$ J(L). Since J(L) $\subseteq$ L$^2$ always holds, L$^2$=J(L).
\vspace{5mm}

Proposition 11 . Let L be a Leibniz algebra over a field of characteristic 0. Let R=R(L). Then J(L)=LR+RL.

\vspace{5mm}

Proof. Let S be a Levi factor of L. Then S is Lie and S=S$_1\bigoplus...\bigoplus$S$_t$ where each S$_i$ is simple. Let M$_i$ be the sum of R and all of the S$_j$ except S$_i$. Then M$_i$ is a maximal ideal of L. Then J(L) $\subseteq \bigcap$ M$_i$=R and J(L) $\subseteq$ R$\cap$ L$^2$.

\vspace{5mm}

If M is a maximal ideal of L, then L/M is abelian or simple. In the first case, L$^2\subseteq$ M and in the second case, R $\subseteq $M. The intersection of all of the first type of M contains L$^2$, and all of the second type contains R. Therefore,  R$\cap$ L$^2 \subseteq$ J(L) and the result follows.

\vspace{5mm}

Corollary 12. J(L) is nilpotent.

\vspace{5mm}

Proof; J(L)=LR+RL $\subseteq$ Nil(L) by Propositions 9 and 11.

\vspace{5mm}

Corollary 13. $\Phi$(L) $\subseteq$ J(L) when L has characteristic 0.

\vspace{5mm}

Proof: L=R+S as in the Levi decomposition and S is Lie. Hence $\Phi$(S)=0. Thus $\Phi$(L) $\subseteq$ R. Since $\Phi$(L) $\subseteq $L$^2$, it follows that $\Phi$(L)$\subseteq$R$\cap$L$^2$=LR+RL=J(L) using Proposition 9 and 11.

\vspace{5mm}

Proposition 14. Let B be a nilpotent ideal in a Leibniz algebra, L. Then J(B) $\subseteq$ J(L).

\vspace{5mm}

Proof: Since B is nilpotent, J(B)=B$^2$, which is an ideal in L. Suppose that x $\in$ J(B), x $\notin$ J(L). Let M be a maximal ideal of L such that x is not in M. Then L=J(B)+M and B=J(B)+(M$\cap$B). M$\cap$B is a proper ideal of B that supplements J(B), a contradiction.

\vspace{5mm}

V. LEIBNIZ ALGEBRAS WITH A UNIQUE MAXIMAL IDEAL.

\vspace{5mm}

Lie algebras with a unique maximal ideal were classified in [5] when the field of scalars has characteristic 0. We extend this result to Leibniz algebras. 

\vspace{5mm}

Lemma 15. Suppose that L=N+$<$x$>$ where $<$x$>$ is a one dimensional vector space, x$^2 \neq 0$, and N is an abelian ideal. Then xL is an ideal in L.

\vspace{5mm}

Proof.  Since xL $\subseteq$ N, N(xL)=(xL)N=0. (xL)x=(x(x+N))x=(x$^2$+xN)x=(xN)x=x(Nx)-N(x$^2$) $\subseteq$ xL since x$^2 \in$ N and N is abelian. Furthermore x(xL) $\subseteq$ xL. Hence the result holds.

\vspace{5mm}

Lemma 16.  Suppose that L=N+$<$x$>$ where $<$x$>$ is a one dimensional vector space, x$^2 \neq 0$ and N is an abelian ideal. If xL+Lx=N, then xL=N.

\vspace{5mm}

Proof. If xL $\neq$ N, then xL is contained in an ideal M of L , M $\subseteq $ N, and N/M is a minimal ideal of L/M. Note that x$^2\in$ M. We may take M=0. Then xN=0, which yields Nx=0 by [1]. Then Lx=(N+x)x =0. Therefore xL+Lx=0, a contradiction

\vspace{5mm}

We now show

\vspace{5mm}

Theorem 17. Let L be a Leibniz algebra over a field of characteristic 0.  L has a unique maximal ideal if and only if one of the following holds.

\vspace{5mm}

1) L is nilpotent and cyclic, or
\vspace{5mm}

2) L=Nil(L)+S where S is simple and N/N$^2$=S(N/N$^2)$+(N/N$^2$)S, or
\vspace{5mm}

3) L=Nil(L)+$<x>$ where $<x>$ is a one dimensional vector space, x$^2 \neq 0$ and N/(N$^2$)=x(L/N$^2$), or
\vspace{5mm}

4) L=Nil(L)+$<x>$  where x$^2$=0 and N/(N$^2$)=x(N/N$^2$)+(N/N$^2$)x.

\vspace{5mm}

Proof.  Let L have only one maximal ideal which is then J(L)=J. Then L/J is simple or one-dimensional. Since J $\subseteq$ Nil(L)=N, if L/J is simple, then J=N and L=J+S=N+S. To show the second part of (2), assume that N$^2$=0 and let R=Rad(L). Then N=R and J=LR+RL=SN+NS.

\vspace{5mm}

Suppose that dim(L/J)=1. Since J $\subseteq$ N $\subseteq$ L, either J=N or N=L.  If N=L, then L is nilpotent and J=$\Phi$(L)=L$^2$, so L is cyclic generated by some a and L is nilpotent and (1) holds.

\vspace{5mm}

Suppose that N=J. Then L=N+$<x>$ where $<x>$ is a one dimensional subspace. Suppose that x$^2 \neq 0$ . By Theorem 3.1 of [4], L/N$^2$ is not nilpotent, hence Nil(L/N$^2$) =N/N$^2$ =J/N$^2$=J(L/N$^2$)=x(L/N$^2$)+(L/N$^2$)x. Suppose that N$^2$=0 and T is the algebra generated by x. Then Nil(L)=J(L)=xL+Lx=T$^2$+xN+Nx. Then J(L/T$^2$)=x(N/T$^2$)+(N/T$^2$)x=x(N/T$^2$). by Lemma 16. Hence Nil(L)= J(L)=T$^2$+xN=xL and (3) holds. Assume that x$^2$=0 and N$^2$=0. As in the last case, Nil(L)=J(L)=LR+RL=xN+Nx.  Hence (4) holds.

\vspace{5mm}

We now show that each of the algebras in 1-4 have a unique maximal ideal. For the algebras in (3), we have the following lemma.

\vspace{5mm}

Lemma 18. Suppose that L is Leibniz, and L=Nil(L)+$<x>$ where x$^2 \neq 0$. Suppose that xL=N and N$^2$=0 where N=Nil(L). Then Nil(L)=J(L). The algebras in (3) have a unique maximal ideal.

\vspace{5mm}

Proof.  Let J=J(L). Since xL=N, N=L$^2$. Let T be the subalgebra generated by x. Since dim(Ker L$_x|_L$)=1 and 0 $\neq$ Ker L$_x|_T \subseteq$ Ker L$_x|_L$, it follows that  Ker L$_x|_T $= Ker L$_x|_L$. Since dim(Ker L$_x|_T$)=1, the Fitting null component of L$_x|_ T$ is the union T$_0 \subseteq$ T$_1\subseteq...$T$_k$ where dim(T$_j$)=j and xT$_j \subseteq $ T$_{j-1}.$ The same holds for the Fitting null component of L$_x$ on L where the invariant subspaces end with L$_m$. Clealy T$_i$=L$_i$ for i $\leq$ k. We claim that m=k. There is a T$_i$ that is not contained in N since the Fitting one component of L$_x$ on T is contained in N. Then $\alpha$x+n $\in $T$_i$, $\alpha \neq 0$ n $\in$ N. If $i\neq$ k, then there exists $\beta$x+p such that x($\beta$x+p)=$\alpha$x+n, which is impossible since x($\beta$x+p) $\in$ N=L$^2$. Thus T$_{k-1} \subseteq$ N while T$_{k}$ is not contained in N. Since T$_k$=L$_k$ is not contained in N, the chain of L$_i$'s must stop at L$_k$ since only the final term in the string is not in N. Hence the Fitting null component for L$_x$ acting on L is the same as L$_x$ acting on T and the Fitting null component of L$_x$ on N and T$^2$ are the same, both equal to T$_{k-1} \subseteq$ T$^2$.

\vspace{5mm}

Let J=J(L). Since xL=N , N=L$^2$. Let N=N$_0$+N$_1$ be the Fitting decomposition of L$_x$ on N. Let K be a maximal ideal of L and suppose that K is not contained in L$^2$. There exists y
$\in$ K such that y=$\alpha$x+n$_0$+n$_1$ where $\alpha \neq 0$, n$_0 \in$ N$_0$ and n$_1 \in$ N$_1$. For any t$\in$ N$_1$, there exists s $\in$ N$_1$ such that xs=t. Then ys=($\alpha$x +n$_0$+n$_1$)s=$\alpha$t since N$^2 = 0$. Thus $\alpha$t $\in$ K and N$_1 \subseteq$ K. Also n$_0 \in$ N$_0$=T$_{k-1} \subseteq$ T$^2$  by the last paragraph. Hence n$_0=\alpha_2$x$^2+...+\alpha_t$x$^t$ and yx=($\alpha$x+$\alpha_2$x$^2+...+\alpha_t$x$^t$ +n$_1$)x =$\alpha$x$^2$ +n$_1$x $\in$ K. Since n$_1$x $\in$K, $\alpha$x$^2 \in $ K. Therefore x$^2 \in$ K and T$^2 \subseteq $ K. Since N$_0 \subseteq $ T$^2$ and N$_1\subseteq$ K, N $\subseteq$ K. Since J(L)=N has codimension 1, J(L)=N=K for all maximal ideals K
of L. Hence algebras as in (3) have a unique maximal ideal.

\vspace{5mm}

We now consider the algebras in (1), (2), and (4). Suppose that L is as in (2). We may assume that N$^2$=0 since N/N$^2$=Nil(L/N$^2$). Thus J=LR+RL=SN+NS=Nil(L) and J(L) is the only maximal ideal of L. If L is as in (4), then assume that N$^2$=0 since N/N$^2$=Nil(L/N$^2$). Then J(L)=LR+RL=xN+Nx. Since dim(L/N)=1, N is the only maximal ideal of L.

\vspace{5mm}

If L is as in (1), then L$^2$=J(L) since L is nilpotent, and J(L) is the unique maximal ideal.

\vspace{10mm}

REFERENCES

\vspace{5mm}

1.  D.W.Barnes, Some theorems on Leibniz algebras, Comm. In Alg. 39. , 2011, 2463-2472 

2. D. W. Barnes, On Levi's theorem for Leibniz algebras, Bull. Aust. Math. Soc. 86, 2012, 184-185

3. C. Batten Ray, L. Bosko-Dunbar, A. Hedges, J.T.Hird, K. Stagg, E. Stitzinger, A Frattini theory for Leibniz algebras, Comm. in Alg. To appear, arXiv 1108.2451

4. C. Batten Ray, A. Combs, N. Gin, A. Hedges, J.T.Hird, L. Zack, Nilpotent Lie and Leibniz algebras, Comm. In Alg. To appear, arXiv 1207.3739

5. P. Benito Clavijo,  Lie algebras in which the lattice of ideals form a chain, Comm. In Alg. 20, 1992, 93-108, 

6. L. Bosko, A. Hedges, J.T. Hird, N. Schwartz, K. Stagg, Jacobson's refinement of Engel's theorem for Leibniz algebras, Involve 4(3), 2011, 293-296

7. N. Jacobson, Lie algebras, Dover, 1979

8. J. Loday, Une version non commutative des algebres de Lie: les algebres de Leibniz, Enseign Math. 39, 1993,  269-293

9. E. I. Marshall, The Frattini Subalgebra of a Lie Algebra, J. Lond. Math. Soc. 42, 1967  416-422

10. K. Stagg and E. Stitzinger, Minimal non-elementary Lie algebras, Proc. Amer. Math. Soc. 137(7), 2010 , 2435-2437

11.  E. Stitzinger, Frattini subalgebras of a class of solvable Lie algebras, Pac. J. Math. 34(1),  1970, 177-182

12. D. A. Towers, A Frattini theory for algebras, Proc. Lond. Math. Soc. 27(3), 1973,  440-462

13. D. A. Towers, Elementary Lie Algebras,  J. Lond. Math. Soc. 7(2),  1973, 295-302

14. D. A. Towers, V. Varea, Elementary Lie Algebras and Lie A-Algebras, J. Alg. 312,  2007, 891-901

\end{document}